\documentclass[conference]{IEEEtran}
\IEEEoverridecommandlockouts
\pdfoutput=1
\usepackage{cite}
\usepackage{amsmath,amssymb,amsthm,epsfig,dsfont,color,subfigure,empheq,graphicx}
\usepackage{enumerate,url,algpseudocode,algorithm,wasysym,epstopdf,verbatim,balance,enumitem}
\usepackage[table]{xcolor}
\DeclareMathOperator{\rank}{rank}

\DeclareMathOperator{\diag}{dg}

\DeclareMathOperator{\vectorize}{vec}

\newtheorem{lemma}{Lemma}

\newtheorem{theorem}{Theorem}

\newcommand \bzero{\mathbf{0}}
\newcommand \bone{\mathbf{1}}
\newcommand \ba{\mathbf{a}}
\newcommand \bb{\mathbf{b}}

\newcommand \be{\mathbf{e}}

\newcommand \bg{\mathbf{g}} 

\newcommand \bp{\mathbf{p}}
\newcommand \bq{\mathbf{q}}
\newcommand \br{\mathbf{r}}

\newcommand \bv{\mathbf{v}}

\newcommand \bx{\mathbf{x}}
\newcommand \by{\mathbf{y}}
\newcommand \bz{\mathbf{z}}
\newcommand \bA{\mathbf{A}}
\newcommand \bB{\mathbf{B}}
\newcommand \bC{\mathbf{C}}

\newcommand \bE{\mathbf{E}}
\newcommand \bF{\mathbf{F}}
\newcommand \bG{\mathbf{G}}
\newcommand \bH{\mathbf{H}}
\newcommand \bI{\mathbf{I}}

\newcommand \bP{\mathbf{P}}

\newcommand \bR{\mathbf{R}}
\newcommand \bS{\mathbf{S}}

\newcommand \bV{\mathbf{V}}

\newcommand \bX{\mathbf{X}}

\newcommand \bZ{\mathbf{Z}}

\newcommand \bdelta{\boldsymbol{\delta}}

\newcommand \bDelta{\mathbf{\Delta}}

\newcommand \mcA{\mathcal{A}}

\newcommand \mcC{\mathcal{C}}

\newcommand \mcG{\mathcal{G}}

\newcommand \mcL{\mathcal{L}}

\newcommand \mcN{\mathcal{N}}

\newcommand \bmcL{\bar{\mathcal{L}}}

\newcommand \tbp{\tilde{\mathbf{p}}}
\newcommand \tbq{\tilde{\mathbf{q}}}

\newcommand \tbv{\tilde{\mathbf{v}}}

\newcommand \bba{\bar{\mathbf{a}}}

\newcommand \bbp{\bar{\mathbf{p}}}

\newcommand \bbA{\bar{\mathbf{A}}}

\newcommand \bbH{\bar{\mathbf{H}}}

\allowdisplaybreaks

\begin{document}
	\title{An MILP Approach for Distribution Grid\\
	Topology Identification using Inverter Probing}
	
	\author{
	\thanks{Work partially supported by the US National Science Foundation CAREER grant 1751085.}
		\IEEEauthorblockN{Sina Taheri, Vassilis Kekatos, and Guido Cavraro}
		\IEEEauthorblockA{Bradley Dept. of ECE, Virginia Tech\\
			Blacksburg, VA 24060\\
			Emails: \{sinataheri,kekatos,cavraro\}@vt.edu}
		}
\maketitle

\begin{abstract}
Although knowing the feeder topology and line impedances is a prerequisite for solving any grid optimization task, utilities oftentimes have limited or outdated information on their electric network assets. Given the rampant integration of smart inverters, we have previously advocated perturbing their power injections to unveil the underlying grid topology using the induced voltage responses. Under an approximate grid model, the perturbed power injections and the collected voltage deviations obey a linear regression setup, where the unknown is the vector of line resistances. Building on this model, topology processing can be performed in two steps. Given a candidate radial topology, the line resistances can be estimated via a least-squares (LS) fit on the probing data. The topology attaining the best fit can be then selected. To avoid evaluating the exponentially many candidate topologies, this two-step approach is uniquely formulated as a mixed-integer linear program (MILP) using the McCormick relaxation. If the recovered topology is not radial, a second, computationally more demanding MILP confines the search only within radial topologies. Numerical tests explain how topology recovery depends on the noise level and probing duration, and demonstrate that the first simpler MILP yields a tree topology in 90\% of the cases tested.
\end{abstract}
\vspace{6pt}
	
\allowdisplaybreaks
	
\section{Introduction}\label{sec:intro}

Distribution grid operators are currently challenged by inaccurate knowledge of the underlying electrical topology. Nonetheless, knowing the network topology is prerequisite for accomplishing any feeder-level optimization task. Some utilities know only their primary network infrastructure and line types, but may not know the energized lines and their precise impedances. Other utilities may have more detailed feeder information, yet it may be outdated. Following a natural disaster, line crew members oftentimes restore grid segments without logging topological changes. This signifies the need for dynamic identification of the feeder topologies.

Several works build on the properties of second-order statistics from smart meter data to infer feeder topologies \cite{Deka18}, \cite{BoSch13}, \cite{ParkDekaChertkov}. A Wiener filtering approach using wide-sense stationary processes on radial networks is put forth in~\cite{TalDeMate2017}. Nonetheless, sample statistics converge to their ensemble values only after a large number of grid data has been collected, thus rendering topology estimates possibly obsolete. 

Detecting which lines are energized can be posed as a maximum likelihood detection task~\cite{VGS17}, \cite{Sharon12}. Granted power readings at all leaf buses and selected lines, topology identification has been formulated as a spanning tree recovery using the notion of graph cycles~\cite{sevlian2015distribution}, while line impedances are estimated via a total least-squares fit in~\cite{patopa}. In \cite{zhao2017learning}, deep neural networks are trained to detect the status of transmission line statuses; nevertheless, the standard PQ/PV power flow dataset feeding the classifiers may not be available in distribution grids. Exploiting the linear relationship between nodal voltage and current phasors, a Kron-reduced admittance matrix is recovered via a low rank-plus-sparse decomposition in~\cite{Ardakanian17}, though the deployment of micro-phasor measurement units occurs at a slower pace in distribution grids. 

The existing schemes rely on passively collected data from smart meters and grid sensors to identify the grid topology. Taking a different approach, we have recently proposed an active data acquisition scheme for load~\cite{BKV17}; or topology recovery~\cite{CaKe18a}, \cite{CaKe18b}. The idea is to intentionally but momentarily perturb the injections of smart inverters and possibly infer the grid topology from the recorded voltages. Perturbing the inverter control loops to infer DC microgrids has been suggested in~\cite{Scaglione2017}. Line impedances have been estimated by having inverters injecting harmonics in \cite{Ciobotaru07}. 

This work improves on \cite{CaKe18a} and \cite{CaKe18b} as follows: Reference \cite{CaKe18a} developed graph algorithms for topology identification presuming that the voltages collected upon inverter probing are noiseless. The devised graph algorithms are applicable to noisy data only for prolonged probing periods. To deal with noisy data in a more practical manner, a convex relaxation approach was proposed in \cite{CaKe18b}, yet without performance guarantees. Aiming at the sweet spot between the two previous topology processing schemes, this work poses grid topology learning as a mixed-integer linear program (MILP) using the powerful technique of McCormick relaxation. If the recovered topology is not radial, we confine the search only within radial topologies using a second, computationally more challenging MILP. Numerical tests explain how topology recovery depends on the level of metering noise and the duration of probing. They also demonstrate that the first MILP yields a tree topology in 90\% of the tested cases and at a shorter computational time.

Regarding \emph{notation}, lower- (upper-) case boldface letters denote column vectors (matrices). Calligraphic symbols are reserved for sets. Symbol $^{\top}$ stands for transposition. Vectors $\mathbf{0}$ and $\mathbf{1}$ are the all-zero and all-one vectors, while $\be_m$ is the $m$-th canonical vector. Symbol $\|\mathbf{x}\|_2$ denotes the $\ell_2$-norm of $\mathbf{x}$ and $\diag(\mathbf{x})$ defines a diagonal matrix having $\mathbf{x}$ on its diagonal.

\section{Problem Formulation}\label{sec:problem}

This section models the data collected via inverter probing after reviewing an approximate grid model.

\subsection{Power Grid Modeling}\label{subsec:model}
A radial single-phase power distribution grid having $N+1$ buses can be modeled by a tree graph $\mcG=(\mcN,\mcL)$. Grid buses are represented by the nodes in set $\mcN:=\{0,\ldots,N\}$, and distribution lines are captured by the edges in $\mcL$. The complex power injection and the voltage magnitude at bus $n$ are denoted by $p_n+jq_n$ and $v_n$, respectively. The substation is indexed by $n=0$ and its voltage is fixed at $v_0=1$. Graph $\mcG$ is rooted at the substation and $|\mcL|=N$. The voltages and power injections at all buses excluding the substation are stacked accordingly in vectors $\bv_t$, $\bp_t$, and $\bq_t$ for time $t$. 

The grid connectivity is captured by the branch-bus incidence matrix. This matrix can be partitioned as $[\ba_0~\bA]$, where its first column $\ba_0$ corresponds to the substation and the rest of its columns form the \emph{reduced incidence matrix} $\bA \in\{0,\pm1\}^{N\times N}$. Matrix $\bA$ is square and invertible for radial grids~\cite{GodsilRoyle}.

Voltages are non-linearly related to power injections. However, by linearizing the power flow equations around the flat voltage profile $\mathbf{1}+j\mathbf{0}$, we obtain the approximate model~\cite{Saverio2}
\begin{equation}\label{eq:model}
\bv_t \simeq \bG^{-1}\bp_t + \bB^{-1}\bq_t + \bone
\end{equation}
where $\bG:=\bA^\top \diag^{-1}(\br)\bA$ and $\bB:=\bA^\top \diag^{-1}(\bx)\bA$; and vector $\br+j\bx$ collects all line impedances. Applying \eqref{eq:model} over two consecutive time instances $t$ and $t-1$ and taking the difference gives
\begin{equation}\label{eq:dmodel}
\tbv_t \simeq \bG^{-1}\tbp_t + \bB^{-1}\tbq_t
\end{equation}
for the differential voltages $\tbv_t:=\bv_t-\bv_{t-1}$; and differential injections $\tbp_t:=\bp_{t}-\bp_{t-1}$ and $\tbq_{t}:=\bq_{t}-\bq_{t-1}$. Building on \eqref{eq:dmodel}, we next elaborate on grid probing.

\subsection{Grid Probing via Inverter Perturbations}\label{subsec:probing}
Let $\mcC \subseteq \mcN$ be the subset of buses hosting smart inverters with cardinality $C:=|\mcC|$. An inverter can be commanded to shed solar generation; (dis)-charge an energy storage unit; or change its power factor. The intentional perturbation of power injections lasts for a second or two, and occurs (a)-synchronously across inverters. 

The perturbation from time $t-1$ to $t$ comprises the $t$-th probing slot. During this slot, one or more inverters are commanded to change their active power while keeping their reactive power unchanged. Let $\delta_{c,t}$ be the change in active power from time $t-1$ to $t$ at bus $c\in\mcC$. Stack all inverter perturbations in vector $\bdelta_t\in\mathbb{R}^C$. Ignoring for now any possible load variations during this 1-sec interval, we get that $\tbq_t=\bzero$ and hence the voltage differences over probing slot $t$ become
\begin{equation}\label{eq:dv}
\tbv_t \simeq \bG^{-1}\bI_\mcC \bdelta_t
\end{equation}
where the $N\times C$ matrix $\bI_\mcC$ collects the columns of the $N\times N$ identity matrix corresponding to the buses in $\mcC$. 

The grid is probed over $T$ probing slots. Stacking the probing slots $\{\bdelta_t\}_{t=1}^T$ as columns of matrix $\bDelta$, and the measured voltage differences $\{\tbv_t\}_{t=1}^T$ as columns of $\bV$ gives
\begin{equation}\label{eq:dV0}
\bV \simeq \bG^{-1} \bI_{\mcC} \bDelta.
\end{equation}
Premultiplying \eqref{eq:dV0} by $\bG$ provides the probing data model
\begin{equation}\label{eq:dV}
\bP=\bG\bV+\bE.
\end{equation}
where $\bP:=\bI_{\mcC} \bDelta$. Matrix $\bE$ captures the error introduced by the approximate grid model; measurement noise; and unmodeled load variations.

If the grid is probed at all buses hosting inverters and voltages are recorded at all buses, could one recover the resistive topology of the grid, that is matrix $\bG$? A sufficient condition ensuring identifiability is repeated here for completeness.

\begin{theorem}[\cite{CaKe18b}]\label{th:leaf+id}
Given noiseless probing data $(\bV,\bDelta)$ with $\bV = \bG \bI_{\mcC} \bDelta$ and $\textrm{rank}(\bDelta)=C$, matrix $\bG$ is identifiable if the grid is probed at all leaf buses.
\end{theorem}

We need $T\geq C$ probing slots for $\bDelta$ to be full row-rank. For example, a diagonal $\bDelta=\bDelta_d$ with $\bDelta_d := \diag\{\delta_m\}$ for non-zero $\delta_m$'s can be implemented asynchronously for $T=C$. Alternatively, if inverters have to be brought back to their original setpoints, one can implement the perturbations $\bDelta=[\bDelta_d~-\bDelta_d]$ within $T=2C$ slots. These two examples showcase that probing can be either operator-instructed for the purpose of topology processing, or occur naturally while the inverters adjust their setpoints.

If not all leaf buses are probed, we can still correctly identify a subgraph of $\mcG$; see~\cite{CaKe18b}: One has to remove the descendant buses of all probed buses to convert $\mcG$ to $\mcG'$. If all the leaf buses of $\mcG'$ are probed, then $\mcG'$ is identifiable. The previous claims hold if voltages are collected at all buses. If voltages are collected only at probing buses, one can still find a reduced graph bearing similarities to the actual grid~\cite{CaKe18a}. We henceforth assume that voltages are collected at all buses. Inverters change their active injections to recover the connectivity and line resistances based on \eqref{eq:dV}. Once the connectivity has been found, line reactances can be readily estimated via reactive probing. The roles of (re)active probing can be apparently interchanged.

The key question is how to use \eqref{eq:dV} to find the grid topology. To this end, reference \cite{CaKe18a} develops graph algorithms that recover the topology exactly under noiseless data, but require long probing periods for noisy data. To cope with the practical scenario of noisy data, reference \cite{CaKe18b} puts forth a convex relaxation, which however lacks performance guarantees. To unveil grid topologies from noisy data within reasonable probing intervals and with performance guarantees, we next devise an MILP formulation.

\section{Topology Processing with Probing Data}\label{sec:model}
To learn the network from \eqref{eq:dV}, we proceed in two steps: \emph{i)} Assuming a topology (matrix $\bA$), we find a least-squares (LS) estimate for the vector of inverse resistances $\bg$; and \emph{ii)} Recover the topology attaining the smallest LS fit on the probing data.

\subsection{Estimating Line Resistances}\label{subsec:resistances}
Suppose that the probing data of \eqref{eq:dV} have been generated from the topology captured by $\bA$. Then \eqref{eq:dV} is a linear measurement model over the inverse resistances, which we denote simply as $\diag(\bg):=\diag^{-1}(\br)$. It is convenient to vectorize the $N\times T$ matrix equation of \eqref{eq:dV} using the $\vectorize(\cdot)$ operator. When this operator is applied on matrix $\bP\in\mathbb{R}^{N\times T}$, it stacks the columns of $\bP$ and returns the vector $\vectorize(\bP)\in\mathbb{R}^{NT}$. A key property from linear algebra is that~\cite{Khatri}
\[\vectorize\left(\bX\diag(\by)\bZ^\top\right)=(\bZ\ast \bX)\by\]
where $\ast$ is the Khatri-Rao matrix product, a column-wise version of the Kronecker product; see \cite{Khatri} for details. 

Thanks to this property, the probing data of \eqref{eq:dV} can be compactly expressed as
\begin{equation}\label{eq:vectorized}
\bp= \bH\bg +\be
\end{equation}
where $\bp:=\vectorize(\bP)$, $\be:=\vectorize(\bE)$, and
\begin{equation}\label{eq:H}
\bH:=\bV^\top\bA^\top\ast\bA^\top.
\end{equation}
Given $\bp$ and assuming $\bH$ to be known, the vector $\bg$ can be estimated using the LS fit
\begin{equation}\label{eq:LS}
\min_{\bg}\left\|\bp-\bH\bg\right\|_2^2.
\end{equation}
If $\bH$ is full column-rank, problem \eqref{eq:LS} has a unique minimizer $\bg^*:=\left(\bH^\top\bH\right)^{-1}\bH^\top\bp$. Lemma~\ref{le:Hrank} ensures that under conditions identical to those of Theorem~\ref{th:leaf+id}, matrix $\bH$ is full column-rank; see the appendix for a proof.

\begin{lemma}\label{le:Hrank}
For any matrix $\bA$ corresponding to a tree topology and voltage data $\bV$ collected by probing the leaf nodes of a radial grid with $\rank(\bDelta)=C$, matrix $\bH$ defined in \eqref{eq:H} is full column-rank.
\end{lemma}

Under the conditions of Lemma~\ref{le:Hrank} and given $\bA$, the minimizer $\bg^*$ of \eqref{eq:LS} attains the LS cost:
\begin{equation}\label{eq:LScost}
\left\|\bp-\bH\bg^*\right\|_2^2 = \|\bp\|_2^2-\bp^\top\bH\left(\bH^\top\bH\right)^{-1}\bH^\top\bp.
\end{equation}
The estimated inverse resistances $\bg^*$ correspond to the distribution lines assumed energized given $\bA$; recall that every row of $\bA$ is associated with a line. 

\subsection{Learning the Grid Topology}\label{subsec:topology}
If the energized topology is unknown, both $\bA$ and $\bg$ need to be found. The grid topology can be found so that the LS cost of \eqref{eq:LScost} is minimized over $\bA$. Under Theorem~\ref{th:leaf+id} and noiseless data, the minimizing $\bA$ agrees with the actual grid topology that gave rise to the voltage data $\bV$. Since $\|\bp\|_2^2$ is fixed, matrix $\bA$ can be recovered by minimizing
\begin{equation}\label{eq:LScostA}
f(\bA):=-\bp^\top\bH\left(\bH^\top\bH\right)^{-1}\bH^\top\bp
\end{equation}
where the dependence on $\bA$ is through $\bH$ as defined in \eqref{eq:H}. 

The operator would like to identify the energized lines from a set $\bmcL$ of \emph{candidate lines} with $\bar{L}:=|\bmcL|$. If the line infrastructure is known, the set $\bmcL$ consists of all actual lines and $\bar{L}$ is typically in the order of $N$. Otherwise, the set $\bmcL$ may include the $\bar{L}=N(N+1)/2$ possible connections among all $N$ buses. One can perform an exhaustive search over all possible spanning trees, which are in general exponentially many. Aiming towards a more systematic search, we next pose the minimization of $f(\bA)$ as a mixed-integer program.

Each candidate line is associated with a row of the \emph{augmented} incidence matrix $[\bba_o~\bbA]\in\{0,\pm1\}^{\bar{L}\times (N+1)}$. Then, the rows of $\bA$ are a subset of the rows of $\bbA$. This row sampling is captured by a selection matrix $\bS\in \{0,1\}^{N\times \bar{L}}$
\begin{equation}\label{eq:bA}
\bA = \bS\bbA.
\end{equation}
Each row of $\bS$ has a single entry equal to one. It is not hard to verify that $\bS$ satisfies the properties
\begin{subequations}\label{eq:bS}
\begin{align}
\bS\bS^\top &= \bI_N\label{eq:bS:1}\\
\bS^\top\bS &= \diag(\bb)\label{eq:bS:2}
\end{align}
\end{subequations}
where $\bb$ is an $\bar{L}$-length indicator vector with $b_\ell=1$, if line $\ell\in\bmcL$ is energized; and $b_\ell=0$, otherwise. Apparently finding $\bA$ amounts to finding $\bS$ or $\bb$. We next express matrix $\bH$ in terms of $\bS$. 

\begin{lemma}\label{le:bH}
For $\bA = \bS\bbA$, the matrix $\bH$ defined in \eqref{eq:H} can be written as
\begin{equation}\label{eq:bH}
\bH = \bbH\bS^\top
\end{equation}
where $\bbH:=\bV^\top\bbA^\top\ast\bbA^\top$.
\end{lemma}

\begin{IEEEproof}
Plugging $\bA = \bS\bbA$ into the definition of $\bH$ yields
\[\bH:=\bV^\top\bA^\top\ast\bA^\top
=\left(\bV^\top\bbA^\top\bS^\top\right)\ast\left(\bbA^\top\bS^\top\right).\]
Because pre-multiplying a matrix by $\bS$ performs row selection, post-multiplying a matrix by $\bS^\top$ performs column selection. Matrix $\bH$ is the Khatri-Rao product of $\bV^\top\bbA^\top$ and $\bbA^\top$, after both being column-sampled by $\bS^\top$. Since the Khatri-Rao applies column-wise, the claim follows. 
\end{IEEEproof}

Based on \eqref{eq:bH}, we next reformulate the matrix inverse of \eqref{eq:LScost}. To this end, pick an $\alpha>0$ and express $\bH^\top\bH$ as
\begin{align}\label{eq:MIL1}
\bH^\top\bH&=\bS\bbH^\top\bbH\bS^\top+\alpha\bI_N-\alpha\bI_N\nonumber\\
&= \bS\left(\bbH^\top\bbH+\alpha\bI_N\right)\bS^\top-\alpha\bI_N\nonumber\\
&=\alpha \bS\bC^{-1}\bS^\top-\alpha\bI_N
\end{align}
where the second equality stems from \eqref{eq:bS:1} and $\bC^{-1}:=\alpha^{-1}\bbH^\top\bbH+\bI_N$ is a given matrix. Thanks to \eqref{eq:MIL1} and property \eqref{eq:bS:2}, the inverse of $\bH^\top\bH$ is amenable to the matrix inversion lemma to provide
\begin{equation}\label{eq:MIL2}
\left(\bH^\top\bH\right)^{-1} =-\alpha^{-1}\bI_N-\alpha^{-1}\bS\left[\bC-\diag(\bb)\right]^{-1}\bS^\top.
\end{equation}

Using \eqref{eq:MIL2} and \eqref{eq:bS:2} again, we can express \eqref{eq:LScostA} as
\begin{equation}\label{eq:alpha_f}
\alpha f(\bA) =\bbp^\top\big[\diag(\bb)+\diag(\bb)\left(\bC-\diag(\bb)\right)^{-1}\diag(\bb)\big]\bbp
\end{equation}
where $\bbp:=\bbH^\top\bp$ is a known vector. Therefore, minimizing $f(\bA)$ is equivalent to minimizing the right-hand side of \eqref{eq:alpha_f} over the binary variable $\bb$, which is a mixed-integer non-convex problem. We will reformulate \eqref{eq:alpha_f} into an MILP.

\section{MILP Solvers}\label{sec:solver}
In search of the actual line status vector $\bb$, adding any prior information as constraints can be helpful. Seeking a tree topology, there should be exactly $N$ energized lines. Moreover, every bus should be connected to at least one energized line. These requirements are encoded via the linear constraints
\begin{subequations}\label{cons:b}
	\begin{align}
	&~\bone_{\bar{L}}^\top\bb = N\label{cons:radial}\\
	&~\left|[\bba_o~\bbA]\right|^\top\bb\geq\bone_{N+1}\label{cons:degree}
	\end{align}
\end{subequations}
with the absolute value understood entry-wise. The matrix-vector product in the left-hand side of \eqref{cons:degree} counts the number of energized lines incident to each bus.

Returning to \eqref{eq:alpha_f} and to bypass the matrix inverse, one can instead try solving the problem
\begin{subequations}\label{eq:stage1}
	\begin{align}
	\underset{\bb,\bz}{\min}~&~f'(\bb,\bz):=\bbp^\top\diag(\bb)\bbp + \bbp^\top\diag(\bb)\bz\label{eq:cost_second}\\
	\mathrm{s.to}~&~(\bC-\diag(\bb))\bz = \diag(\bb)\bbp\label{cons:a_equal}\\
	&~\eqref{cons:b}\label{eq:stage1:c}
	\end{align}
\end{subequations}
where $\bz$ is an auxiliary optimization variable. If matrix $\bC-\diag(\bb)$ is invertible at optimality, then $\bz=\left(\bC-\diag(\bb)\right)^{-1}\diag(\bb)\bbp$ from \eqref{cons:a_equal} and the optimal cost of \eqref{eq:stage1} yields the minimal $\alpha f(\bA)$. Otherwise, the solution to \eqref{eq:stage1} does not necessarily minimize \eqref{eq:alpha_f}. 

Nonetheless, due to its relative simplicity, one can still try solving \eqref{eq:stage1}. Problem \eqref{eq:stage1} involves the products $\diag(\bb)\bz$ between binary and continuous variables. By applying a McCormick relaxation on each bilinear term (see e.g., \cite{Bhela19a}), problem \eqref{eq:stage1} can be expressed as an MILP. The McCormick relaxation requires that the continuous variable is known to be lie in a box. Although we were not able to derive a box interval for $\bz$ analytically, the bound $\|\bz\|_\infty\leq\|\bbp\|_\infty$ worked well during the numerical tests of Section~\ref{sec:test}. 

If problem \eqref{eq:stage1} yields a singular matrix $\bC-\diag(\bb)$ at optimality, one should follow a more elaborate approach: Equation~\eqref{eq:MIL2} suggests that matrix $\bC-\diag(\bb)$ is invertible if and only if $\bH$ is full column rank. From Lemma~\ref{le:Hrank}, we know that $\bH$ is full column-rank if the candidate $\bA$ corresponds to a tree. To confine our search to $\bA$'s corresponding to trees: \emph{i)} we include $\bS$ into the optimization variables; \emph{ii)} introduce an additional variable $\bF$; and \emph{iii)} add the constraint
\begin{equation}\label{eq:FA}
-\bF\bS\bbA = \bI_N.
\end{equation}
Constraint \eqref{eq:FA} dictates that $\bA$ is invertible with $\bA^{-1}=-\bF$; recall $\bS\bbA=\bA$ from \eqref{eq:bA}. Interestingly, matrix $\bF$ has binary entries for any $\bA$ related to a tree~\cite{SGC15}.

Having added $\bS$ in the optimization variables, we also append the linear constraints 
\begin{subequations}\label{eq:morecon}
	\begin{align}
	\bS^\top\bone_N &= \bb\label{eq:morecon:1}\\
	\bS\bone_{\bar{L}}&=\bone_N\label{eq:morecon:2}
	\end{align}
\end{subequations}
which surrogate \eqref{eq:bS} and relate variables $\bS$ with $\bb$. 

Putting the pieces together, minimizing \eqref{eq:alpha_f} is equivalent to solving the problem
\begin{align}\label{eq:stage2}
\underset{\bb,\bz,\bS,\bF}{\min}~&~f'(\bb,\bz)\\
\mathrm{s.to}~&~\eqref{cons:b},\eqref{cons:a_equal},\eqref{eq:FA},\eqref{eq:morecon}.\nonumber
\end{align}
Similarly to \eqref{eq:stage1}, the product $\bF\bS$ in \eqref{eq:FA} can be converted to linear inequality constraints through McCormick relaxation. Therefore, problem \eqref{eq:stage2} can be also posed as an MILP.

Constraint~\eqref{eq:FA} becomes computationally expensive for increasing $\bar{L}$, since the number of bilinear terms grows as $\bar{L}N^3$ and the corresponding McCormick constraints grow as $4\bar{L}N^3$. It is exactly for this reason that we suggest solving \eqref{eq:stage2} only if \eqref{eq:stage1} fails to yield an invertible $\bC-\diag(\bb)$ at optimality. In fact, problem \eqref{eq:stage1} returned a tree topology for more than $90\%$ of the numerical tests described below.

\section{Numerical tests}\label{sec:test}
Our topology recovery schemes of \eqref{eq:stage1} and \eqref{eq:stage2} were tested on the IEEE 13-bus system, modified to represent a single-phase network. One may draw $N(N+1)/2=78$ connections between $N=12$ buses and the substation. For each test, the set $\bmcL$ was constructed by appending $12$ extra lines to the $12$ existing lines, so that $\bar{L}=24$. Complying with Theorem~\ref{th:leaf+id}, every leaf node is assumed to host a smart inverter sized for the full load of the respective bus.

Each probing action lasted for a second. To model load variations, a zero-mean normal random variable with standard deviation of $\sigma_\ell = 0.0068$~pu was considered~\cite{Guido_Arghandeh18}. Instead of the approximate model of \eqref{eq:model}, the voltages induced by probing were calculated using the MATPOWER AC power flow solver~\cite{MATPOWER}. To account for metering noise, each voltage magnitude reading $v_n$ was corrupted by multiplicative noise as $v_n(1+\epsilon_n)$, where $\epsilon_n$ was a zero-mean Gaussian with a standard deviation of $\sigma_\epsilon$. Smart meters exhibit an accuracy of $0.2-0.5\%$ for voltage magnitudes, while the accuracy for micro-phasor measurements units is $0.01\%$. For this reason, the parameter $\sigma_\epsilon$ was tested within the range of $5\cdot10^{-5}$ to $10^{-3}$~pu. 

Only one inverter is probed in each probing action to minimize the impact of probings. Each inverter is probed by first nulling its injection and then bringing it back on in the next probing action. For this reason, each inverter is probed even number of times, i.e. $T_\text{inv}$, and we have $T=CT_\text{inv}$ where $C=6$ for our testbed. Because the differential voltages of \eqref{eq:dv} laid in the range of $10^{-3}$, the voltage data were scaled by $1,000$, while powers were not scaled assuming that the inverse resistances were reversely scaled, see~\eqref{eq:dV}. This was only to numerically condition the data fed to the MILP solver of Gurobi~\cite{gurobi}. Using tighter bounds on the continuous variables in McCormick relaxations can accelerate significantly the timing of MILP solvers. Through our experiments, we numerically observed that $\|\bz\|_\infty\leq\|\bbp\|_\infty$ when $\alpha\leq 25$.

\begin{figure}[htbp]
	\centering
	\includegraphics[width=0.4\textwidth]{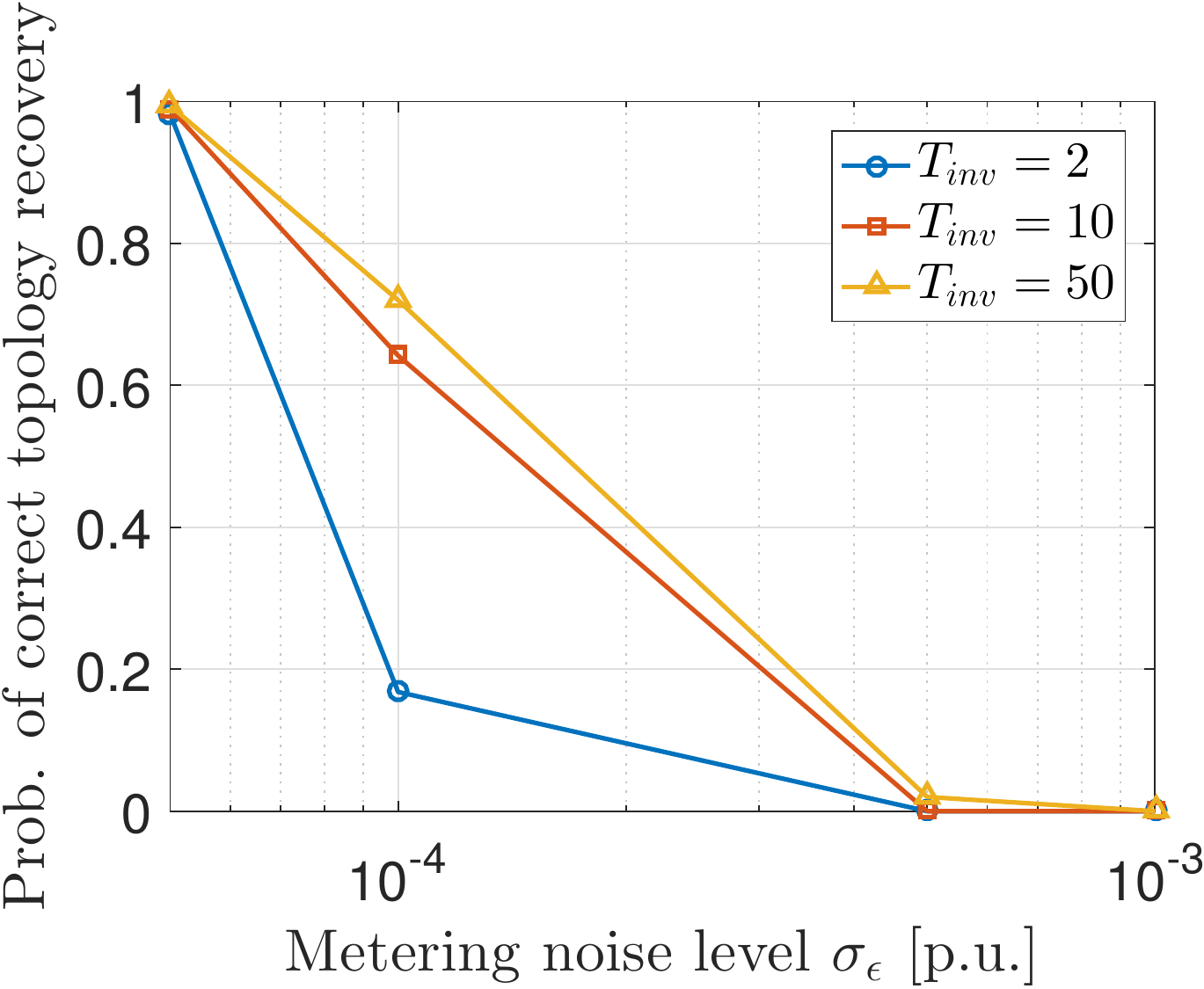}
	\caption{Probability of recovering the correct topology for different noise levels $\sigma_\epsilon$ and number of probing actions per inverter $T_{inv}$.}
	\label{fig:prob_suc}
\end{figure}

We first evaluated the probability of finding the true topology. Figure~\ref{fig:prob_suc} shows the results obtained for different noise levels and number of probing actions per inverter $T_\text{inv}$. The curves indicate that for relatively small $\sigma_\epsilon$, having more probing actions does not improve the probability of topology recovery significantly, though longer probing becomes more effective at higher noise levels. Since misplacing even a single line is counted as failure, this probability metric goes close to zero for higher $\sigma_\epsilon$'s.

\begin{figure}[t]
	\centering
	\includegraphics[width=0.4\textwidth]{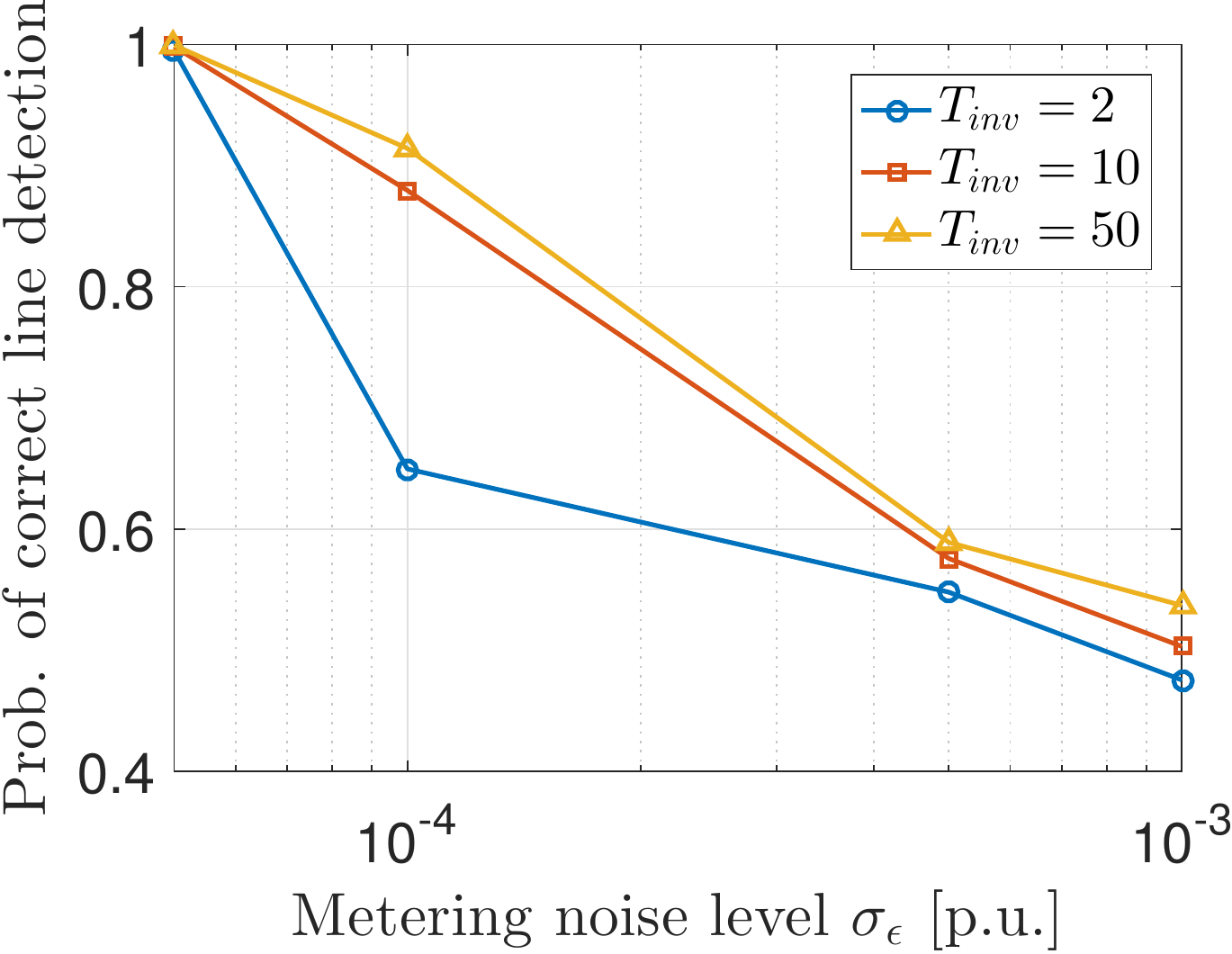}
	\caption{Probability of correct line detection for different noise levels $\sigma_\epsilon$ and number of probing actions per inverter $T_{inv}$.}
	\label{fig:avg_line}
\end{figure}

We also evaluated the probability of correct line detection defined as $\|\bb^*-\bb_o\|_1/\bar{L}$, where $\bb_o$ ($\bb^*$) is the actual (recovered) line status vectors. As demonstrated by Figure~\ref{fig:avg_line}, although the grid topology may not be fully recovered for higher noise levels, the majority of the lines are detected correctly. Because the line status errors for noise levels of $\sigma_\epsilon=5\cdot 10^{-5}$ and $10^{-4}$ were very few, the number of Monte Carlo tests for that range was increased to 1000.

\begin{table}[t]
	\renewcommand{\arraystretch}{1.1}
	\caption{Mean Square Error for Line Inverse Resistances [\%]}
	\vspace*{-1em}
	\label{tbl:MSE} \centering
	\begin{tabular}{|r|r|r|r|r|}
		\hline\hline
		Metering level ($\sigma_\epsilon$)  & $T=2$ & $T=10$ & $T=50$\\
		\hline\hline
		$5\times10^{-5}$~pu & $11.47$ & $10.86$ & $8.63$\\ 
		\hline
		$10^{-4}$~pu   & $35.36$ & $32.11$ & $29.41$\\
		\hline\hline
	\end{tabular}
\end{table}

For the cases where energized lines were correctly identified, we also calculated the mean square error of the estimated inverse resistances. Table~\ref{tbl:MSE} demonstrates that prolonged probing intervals improve estimation accuracy as expected.

Finally, we evaluated the running times of~\eqref{eq:stage1} and~\eqref{eq:stage2}. The tests were performed on a personal computer with Intel Core i7 @ 3.4 GHz (16 GB RAM) using MATLAB and the Gurobi solver. The median running times are shown in Table~\ref{tbl:time}. These results show that the MILP of~\eqref{eq:stage1} is much faster than the MILP of~\eqref{eq:stage2}; this justifies our two-stage approach. 

\begin{table}
	\renewcommand{\arraystretch}{1.1}
	\caption{Running Time for MILPs [sec]}
	\label{tbl:time}
	\vspace*{-1em}
	\label{tbl:Comp_time} \centering
	\begin{tabular}{|r|r|r|r|}
		\hline\hline
		\# of candidate lines $\bar{L}$  & $24$ & $36$ & $48$\\
		\hline\hline
		Simpler MILP of \eqref{eq:stage1} & $1$ & $27$ & $200$\\ 
		\hline
		Exact MILP formulation of \eqref{eq:stage2}  & $47$ & $99$ & $1,035$\\
		\hline\hline
	\end{tabular}
\end{table}

\section{Conclusions}
This work builds upon the active data acquisition paradigm of grid probing. An electric is assumed to be probed at all terminal buses and its voltage response is collected at all buses. For every candidate topology, the related resistances have been found via a simple LS fit. Plugging the estimated resistance values, the topology attaining the smallest LS fit for the probing data is deemed as the actual topology. To avoid enumeration of all possible topologies, the topology identification task has been posed as a mixed-integer non-convex program. Using McCormick linearization, the latter program has been reformulated as two MILPs. The first MILP relaxes the search space into all topologies, whereas the latter searches only radial topologies, and is hence equivalent to the original mixed-integer non-convex program. Numerical tests demonstrate correct topology identification for relatively small measurement noise. For higher noise levels, statuses are correctly detected for the majority of the lines. Both detection probabilities improve as probing increase. Extending the approach to multi-phase systems and combining probing with smart meter data are left for future work. Appending additional constraints derived from graph theoretic properties to the MILP formulation is also an interesting direction.

\appendix[Proof of Lemma~\ref{le:Hrank}]
We first review some properties of the symmetric matrix $\bR:=\bG^{-1}$. In a radial grid, every bus $n$ is connected to the substation through a unique sequence of distribution lines, which we term the \emph{feeding path}. The $(m,n)$-th entry of $\bR$ relates to the intersection of the feeding paths for buses $m$ and $n$; see \cite{Deka18}, \cite{CaKe18a}. In fact, the entry $R_{mn}$ sums up the resistances for the lines belonging to the intersection of the two feeding paths. Therefore, it holds that
\begin{equation}\label{eq:R}
R_{nn}\geq R_{nm}\quad \forall n,m\in\mcN
\end{equation}
with strict inequality if bus $n$ is a leaf. The next lemma establishes another key property of $\bR$.

\begin{figure}[t]
	\centering
	\includegraphics[width = 0.18\textwidth]{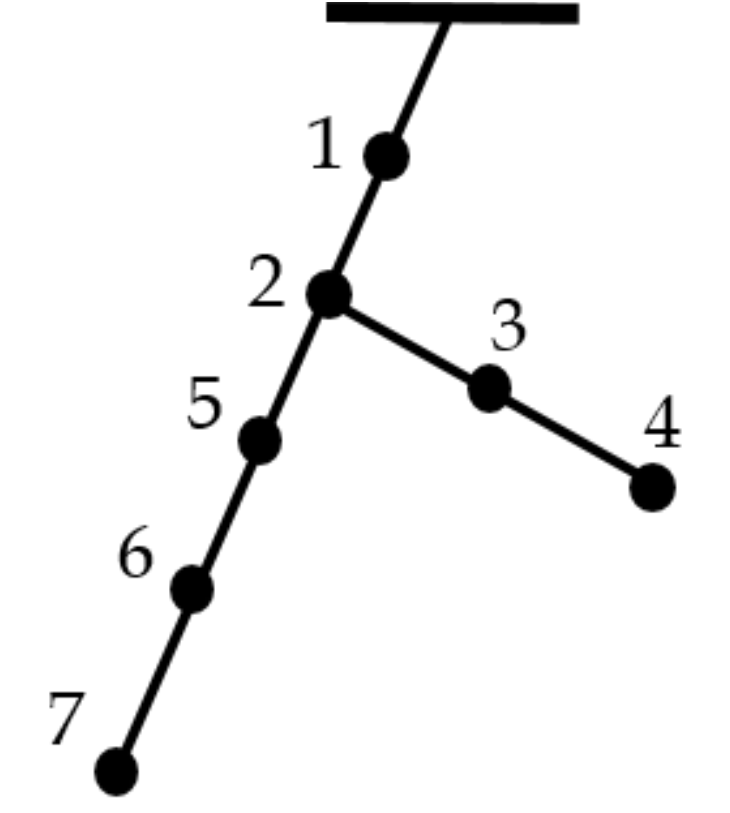}
	\caption{Bus $7$ is a leaf and bus $3$ is not, so $R_{77}>R_{37}$. Buses $2$ and $6$ are not leaves, while $2\in\mcA_6$ and $R_{76}>R_{72}$. Buses $3$ and $6$ are not leaves, but $6\notin\mcA_3$ and $ 3\notin\mcA_6$, whereas $2\in\mcA_6$ and $2\in\mcA_3$, so that $R_{67}>R_{37}$.}
	\label{fig:sample_graph}
\end{figure}

\begin{lemma}\label{le:Rrow}
Matrix $\bR_\mcC :=\bR\bI_\mcC$ has distinct rows if $\mcC$ contains all the leaf nodes of its associated grid.
\end{lemma}

\begin{IEEEproof}
It suffices to show that every pair $(n,m)$ of rows of $\bR_\mcC$ differs by at least one entry. If bus $n$ is a leaf, then $R_{nn}>R_{nm}=R_{mn}$ for all $m\in\mcN$ due to \eqref{eq:R}. Thus, the $n$-th row of $\bR_\mcC$ differs by at least one entry from every other row of $\bR_\mcC$. 

We then have to consider only the row pairs $(n,m)$ between non-leaf buses. Define the ancestors $\mcA_n$ for bus $n\in\mcN$ as the set of buses visited by its feeding path. Based on ancestors, three cases can be identified:
\renewcommand{\labelenumi}{\emph{c\arabic{enumi})}}
\begin{enumerate}
\item If $n\in\mcA_m$, there exists a leaf bus $s\in\mcC$ for which $m\in\mcA_s$, and consequently $ R_{ms}>R_{ns}$.
\item If $m\in\mcA_n$, there exists a leaf bus $s\in\mcC$ for which $n\in\mcA_s$, and consequently $ R_{ns}>R_{ms}$.
\item Otherwise, buses $n$ and $m$ have at least one common ancestor. Let $k$ be the deepest common ancestor of $n$ and $m$. There exists a leaf bus $s\in\mcC$ for which $n\in\mcA_s$ and so $k\in\mcA_s$. Since both $n$ and $k$ are on the feeding path of $s$ and $k\in\mcA_n$, it holds $R_{ns}>R_{ks}$. Bus $k$ is also the deepest common ancestor of $s$ and $m$ and so $R_{ms} = R_{ks}$, which means that $R_{ns}>R_{ms}$.
\end{enumerate}
Figure~\ref{fig:sample_graph} illustrates the different cases.
\end{IEEEproof}

\begin{IEEEproof}[Proof of Lemma~\ref{le:Hrank}]
Matrix $\bH$ is full column-rank if and only if $\bH^\top\bH$ is invertible. From the properties of the Khatri-Rao product~\cite{Khatri}:
\begin{align}\label{eq:bHadamard}
\bH^\top\bH &= \left(\bV^\top\bA^\top\ast\bA^\top\right)^\top\left(\bV^\top\bA^\top\ast\bA^\top\right)\nonumber\\
&=\left(\bA\bV\bV^\top\bA^\top\right)\odot(\bA\bA^\top)
\end{align} 
where $\odot$ denotes the Hadamard (entry-wise) product between two matrices. From~\cite[Th.~5.2.1]{Topics_Horn}, matrix $\bH^\top\bH$ is strictly positive definite and hence invertible, if $\bA\bA^\top\succ \bzero$ and $\bA\bV\bV^\top\bA^\top$ has no zero diagonal entries. The condition $\bA\bA^\top\succ \bzero$ holds since $\bA$ is invertible for radial grids.
		
The $\ell$-th diagonal entry of $\bA\bV\bV^\top\bA^\top$ can be expressed as
\begin{equation}\label{eq:dHadamard}
\left[\bA\bV\bV^\top\bA^\top\right]_{\ell\ell} = \|\ba_\ell^\top\bV\|_2^2
\end{equation}
where $\ba_\ell^\top$ is the $\ell$-th row of $\bA$. In a noiseless probing setup, it holds that $\ba_\ell^\top\bV = \ba_\ell^\top\bG^{-1}\bI_\mcC\bDelta$. Because $\bDelta$ is full row-rank, vector $\ba_\ell^\top\bV$ is zero if and only if the vector $\ba_\ell^\top\bG^{-1}\bI_\mcC$ is zero. Using Lemma~\ref{le:Rrow}, it is evident that $\ba_\ell^\top\bG^{-1}\bI_\mcC= \ba_\ell^\top\bR_\mcC\neq\bzero^\top$, and so matrix $\bA\bV\bV^\top\bA^\top$ in \eqref{eq:bHadamard} has non-zero diagonal entries.
\end{IEEEproof}

\balance
\bibliographystyle{IEEEtran}
\bibliography{myabrv,power}
\end{document}